\documentclass[conference]{IEEEtran}
\IEEEoverridecommandlockouts
\pdfoutput=1
\usepackage{cite}
\makeatletter
\newcommand{\removelatexerror}{\let\@latex@error\@gobble}
\makeatother
\usepackage{hyperref}
\hypersetup{
    colorlinks=true,            
    linkcolor=blue,             
    filecolor=magenta,          
    urlcolor=cyan,              
    pdftitle={Overleaf Example},
    pdfpagemode=FullScreen,
    }
\usepackage{amsmath,amssymb,amsfonts}
\usepackage{graphicx}
\usepackage{subfigure}
\usepackage{textcomp}
\usepackage{xcolor}

\usepackage{algorithmic}
\usepackage{bm}
\usepackage{diagbox}
\usepackage{booktabs}
\usepackage{url}
\usepackage{multirow}
\usepackage{geometry}
\newgeometry{
  top=72pt,
  bottom = 54pt,
  left = 54pt,
  right = 54pt
}
\usepackage{stfloats}
\usepackage[lined,boxed,commentsnumbered,ruled]{algorithm2e}
\def\BibTeX{{\rm B\kern-.05em{\sc i\kern-.025em b}\kern-.08em
    T\kern-.1667em\lower.7ex\hbox{E}\kern-.125emX}}

\title{Enabling Fast Unit Commitment Constraint Screening via Learning Cost Model}

\author{\IEEEauthorblockN{Xuan He\IEEEauthorrefmark{1},
Honglin Wen\IEEEauthorrefmark{2}, Yufan Zhang\IEEEauthorrefmark{2} and
Yize Chen\IEEEauthorrefmark{1}}

    \IEEEauthorblockA{\IEEEauthorrefmark{1}Hong Kong University of Science and Technology (Guangzhou)
    \\xhe085@connect.hkust-gz.edu.cn, yizechen@ust.hk}
    \IEEEauthorblockA{\IEEEauthorrefmark{2}Shanghai Jiao Tong University
    \\\{linlin00, zhangyufan\}@sjtu.edu.cn}}

\begin{document}

\maketitle

\begin{abstract}
    Unit commitment (UC) are essential tools to transmission system operators for finding the most economical and feasible generation schedules and dispatch signals. Constraint screening has been receiving attention as it holds the promise for reducing a number of inactive or redundant constraints in the UC problem, so that the solution process of large scale UC problem can be accelerated by considering the reduced optimization problem. Standard constraint screening approach relies on optimizing over load and generations to find binding line flow constraints, yet the screening is conservative with a large percentage of constraints still reserved for the UC problem. In this paper, we propose a novel machine learning (ML) model to predict the most economical costs given load inputs. Such ML model bridges the cost perspectives of UC decisions to the optimization-based constraint screening model, and can screen out higher proportion of operational constraints. We verify the proposed method's performance on both sample-aware and sample-agnostic setting, and illustrate the proposed scheme can further reduce the computation time on a variety of setup for UC problems. 
\end{abstract}

\section{Introduction}
Solving large-scale optimization problems is one of the cornerstones for many power system operation tasks, such as unit commitment (UC), optimal power flow (OPF), and capacity expansion planning. Many such applications require to solve the problem in a timely manner, as the solutions are essential to market clearing, power grid reliability, and power grid operations~\cite{glover2012power}. And in particular, as the DC power flow is often utilized in UC models, while the generators' ON/OFF statutes are often described as integer variables, UC problems are usually formulated as mixed integer programming (MIP). There is a lot of computation burden to find the most economical UC solutions for such NP-hard problems on large-scale power networks. 

The computational burden of the UC problem can be significant, especially for large-scale transmission systems. This motivates the research on seeking surrogate models that are much smaller than the original UC problems yet ensure equivalence of constrains' binding situation at optimal solution which should be the same as the original one. On the one hand, the active constraints can be kept with the constraint generation technique, where the violated constraints are gradually added into the surrogate model until the solution to the surrogate model is feasible to the original UC problem \cite{xavier2019transmission}. On the other hand, the redundant constraints can be screened out 
 by solving a relaxed optimization for each constraint to identify whether the upper or lower bound of each constraint is redundant \cite{zhai2010fast}. After securely screening constraints, the promise is to reduce the computation time significantly for the reduced UC problems which only involve a subset of original physical constraints in the UC problem. 
 
 Although the total number of UC constraints is theoretically prohibitive with the above methods, empirical evidence shows that a vast majority of the constraints are still redundant and only a smaller subset of constraints could be binding~(equality holds) given the region of load profiles~\cite{ardakani2013identification}. Recent efforts have explored the potential of using cost-driven \cite{porras2021cost} and data-driven \cite{pineda2020data} ways to handle these issues. \cite{porras2021cost} adds an operational cost upper bound to the standard relaxed optimization so as to further narrow the subset of constraints; \cite{pineda2020data} proposes to use the k-nearest neighbor as a screening model to classify the constraint binding statuses using historical samples, yet both methods have tradeoffs in accuracy and efficiency.
 
 Indeed, the availability of historical power system operation records can provide rich information regarding operational decisions, constraint patterns and load profiles \cite{ardakani2018prediction, xavier2021learning}. Note that the system operating pattern corresponds to a unique region of load. The reexamination of the inactive or redundant constraints is needed, once the system operating pattern changes due to the change of the load. As such, \cite{porras2021cost} determines the cost upper bound for different aggregated net demand, i.e., load level. \cite{roald2019implied} considers the constraint screening for varying load ranges and results in a surrogate model that is applicable in a long operation period. In \cite{awadalla2022influence}, the spatial correlation between nodal demands is taken into consideration for identifying the umbrella constraints.

\begin{figure*}[htbp]
\vspace{-1em}
\includegraphics[height=4cm, width=18cm]{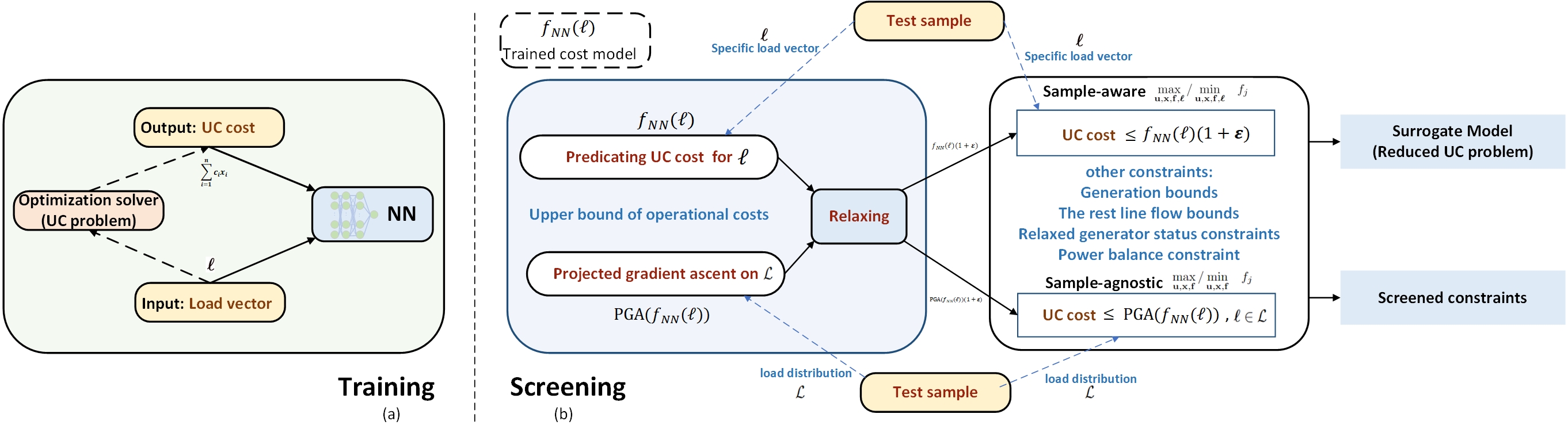}
\centering
\caption{The schematic of our proposed method. We use the historical records involving load vector and UC cost to train a neural network model serving as cost model (a), and further use the trained model to get the upper bound of UC cost so as to screen constraints and obtain a surrogate model (b).} \label{Framework}
\vspace{-1em}
\end{figure*}

However, the data-driven method used in \cite{pineda2020data,ardakani2018prediction, xavier2021learning} directly using machine learning (ML) predictions to classify if constraints are redundant, which may have simplicity but be hard to respond to the change of load region and without guarantee of equivalence. Meanwhile, the cost-driven method proposed in \cite{porras2021cost} can promise the solution accuracy, but it adds integer constraints to the constraint screening problem, making the screening procedure cumbersome to solve. Thus, it is required to come up with a method considering load region to achieve the accuracy-simplicity balance for the constraints screening problem.

In this work, we investigate the potential of designing cost-driven paradigms to efficiently screen constraints in the unit commitment problem. We try to bridge the insights given by machine learning algorithms and the standard optimization-based constraint screening processes. Our approach also utilizes the historical data, but here ML prediction is used to help optimization-based method screen each constraint more efficiently. Specifically, we train a neural network to predict the optimal costs given load inputs~\cite{chen2022learning}. Then we can conveniently upper-bound the search space of constraint screening problem by integrating the cost level predicted and optimized via the trained neural network model, where the formulation of the screening can be treated as a simple linear programming problem. Besides, our method can be flexibly integrated to screen constraints for either one specific load vector or for a given region of load, which we term as \emph{sample-aware}~\cite{zhai2010fast} and \emph{sample-agnostic}~\cite{roald2019implied} constraint screening respectively. In sample-aware case, our proposed method can screen out about 90.03\% of the redundant constraints. In sample-agnostic case, we can realize the UC cost prediction with relative error less than 1\% and remove the redundant constraints without cost error along with saving the solution time.

\vspace{-5pt}
\section{Problem Setup}
\subsection{UC Problem Formulation}

In this paper, we assume the system operators need to decide both the ON/OFF statuses as well as dispatch level for all generators. As the realistic UC problem requires to take start-up and shut-down costs and logic constraints as well as ramp constraints into considerations, which make the analysis of multi-step constraints more complicated, we firstly consider the single-period UC problem as follows:
\begin{subequations}
\label{UC}
    \begin{align}
J(\boldsymbol{\ell})=\min _{\mathbf{u}, \mathbf{x}, \mathbf{f}}\quad  & \sum_{i=1}^n c_i x_i \label{UC:obj}\\
\text { s.t. } \quad & u_i \underline{x}_i \leq x_i \leq u_i \bar{x}_i, \quad i =1, ..., n \label{UC:gen}\\
&-\overline{\mathbf{f}} \leq \mathbf{K f} \leq \overline{\mathbf{f}} \label{UC:flow}\\
& \mathbf{x}+\overline{\mathbf{A}} \mathbf{f}=\boldsymbol{\ell} \label{UC:balance}\\
& u_i \in \{0, 1\}, \quad i=1,...,n \label{UC:u}.
\end{align}
\end{subequations}

In the UC problem, we optimize over the generator statuses $\mathbf{u}$, the generator dispatch $\mathbf{x}$ and the line power flow $\mathbf{f}$ to find the least-cost solutions with cost denoted as $J(\boldsymbol{\ell})$ in the objective function \eqref{UC:obj}. $c_i$ denotes the cost coefficient. Constraint \eqref{UC:gen}, \eqref{UC:flow} and \eqref{UC:balance} denotes the generation bound, the flow bound and the nodal power balance respectively. Note that the power flows are modeled as a DC approximation, while the phase angles are absorbed into the fundamental flows $\mathbf{f}\in \mathbb{R}^{n-1}$~\cite{chen2022learning, bertsimas1997introduction}; $K$ and $\bar{\mathbf{A}}$ map such fundamental flows to flow constraints and nodal power balance respectively. \eqref{UC:u} enforces the binary constraint of generator statuses, where $u_i=1$ indicates the generator is on.

\vspace{-0.5em}
\subsection{Constraint Screening}
Since there are many redundant line flow constraints when seeking the optimal solution of UC problem with given load region or specific load vector, which brings unnecessary computation burden, constraint screening for the line flow constraints can be meaningful. Similar to \cite{zhai2010fast}, we relax the integer variables $\mathbf{u}$ in \eqref{UC} as continuous variables in $[0,1]$, and the screening approach requires to iteratively solve the relaxed optimization problem to find the upper and lower flow values on each transmission line. If the upper and lower bound cannot be reached by the relaxed optimization problem, we can safely screen out that line flow constraint.

For the case that the load region $\mathcal{L}$ is known, a \emph{sample-agnostic constraint screening problem} can be formulated for a group of operating scenarios, which can be given as follows,
\begin{subequations}
\vspace{-1.5em}
\label{screening1}
    \begin{align}
\max_{\mathbf{u}, \mathbf{x}, \mathbf{f}, \boldsymbol{\ell}} / \min _{\mathbf{u}, \mathbf{x}, \mathbf{f}, \boldsymbol{\ell}}\quad  & f_{j} \label{Screening:obj}\\
\text { s.t. } \quad & u_i \underline{x}_i \leq x_i \leq u_i \bar{x}_i, \quad i =1, ..., n \label{Screening:gen}\\
&-\overline{\mathbf{f}}_{\mathcal{F}/j} \leq \mathbf{K}_{\mathcal{F}/j} \tilde{\mathbf{f}} \leq \overline{\mathbf{f}}_{\mathcal{F}/j}\label{Screening:flow}\\
& \mathbf{x}+\overline{\mathbf{A}} \mathbf{f}=\boldsymbol{\ell} \label{Screening:balance}\\
& 0\leq u_i \leq 1, \quad i=1,...,n \label{Screening:u}\\
&  \boldsymbol{\ell} \in \mathcal{L} \label{Screening:load};
\end{align}
\end{subequations}
where $\mathcal{F}/j$ denotes all remaining entries of vectors or matrix which excludes those correspond to $f_j$.

On the contrary, when the specific load vector is available, we can conduct the following \emph{sample-aware constraint screening}:
\begin{subequations}
\label{screening2}
    \begin{align}
\max_{\mathbf{u}, \mathbf{x}, \mathbf{f}} / \min _{\mathbf{u}, \mathbf{x}, \mathbf{f}}\quad  & f_{j}\\
\text { s.t. } \quad & u_i \underline{x}_i \leq x_i \leq u_i \bar{x}_i, \quad i =1, ..., n \label{Screening2:gen}\\
&-\overline{\mathbf{f}}_{\mathcal{F}/j} \leq \mathbf{K}_{\mathcal{F}/j} \tilde{\mathbf{f}} \leq \overline{\mathbf{f}}_{\mathcal{F}/j}\label{Screening2:flow}\\
& \mathbf{x}+\overline{\mathbf{A}} \mathbf{f}=\boldsymbol{\ell} \label{Screening2:balance}\\
& 0\leq u_i \leq 1, \quad i=1,...,n \label{Screening2:u};
\end{align}
\end{subequations}
where $\boldsymbol{\ell}$ is a known load vector for UC problem.

The above formulations are both optimization-based approaches, which seek to find the limit of the flow while keeping all other flow and generation constraints satisfied. However, this approach still allows some line flow values causing unrealistic cost to reach the upper or lower bounds, and thus there are more redundant constraints reserved\cite{porras2021cost}. 

Therefore, it is interesting to consider the economical goal in the original UC problem, minimizing the system cost, to further safely screen out constraints. 

\vspace{0.5em}
\section{Learning to Predict UC Costs}
\subsection{Learning Cost Predictors for Unit Commitment Problem}

As mentioned before, screening without cost objectives enlarge the possible value range of load variables $f_j$, which leads to conservative screening and keep more constraints as non-redundant. To close such gap, in this paper, we investigate if it is possible to tighten the search space of constraint screening  by adding a cost constraint in the form of $\sum_{i=1}^n c_ix_i \leq \Bar{C}$, where $\Bar{C}$ is the upper bound whose value needs to be determined in the following sections.

To achieve such goal, the adopted method should approximate the map between load input and system costs $J(\boldsymbol{\ell})$ well along with predicting system costs efficiently. Thus, in this paper, we use a neural network (NN) to find the upper bound. To train the NN model, we utilize the past record of UC solutions and the training loss between the output of NN model and the actual cost defined as follows,
\begin{align}  \label{Training: loss}
L: =\left\|f_{NN}(\boldsymbol{\ell}) - J(\boldsymbol{\ell})\right\|_{2}^{2};
\end{align}
where $f_{NN}(\boldsymbol{\ell})$ denotes the NN model given load inputs.


In the next subsection, we detail how to connect NN's predicted costs to the constraint screening problems.

\subsection{Tightening the Search Space for Constraint Screening}
Note that the ML model is not directly applied for making operation or dispatch decisions, and alternatively, we are treating the ML prediction as a constraint to reduce the search space of optimization-based constraint screening problems \eqref{screening1} and \eqref{screening2}. With such design, the resulting optimization problem can still find feasible decisions for the original UC problem. We can then add the neural network's prediction as an additional constraint to the original constraint screening problem to further restrict the search space for each transmission's flow bounds. 

In sample-agnostic case, to ensure feasibility of the constraint screening problem after adding the cost constraint for the whole load space along with restricting the searching space effectively, we need to find a predicted cost given by NN which can serve as the upper bound. Then projected gradient ascent~(PGA) algorithm can be adopted to achieve this goal. PGA can find the upper bound iteratively by moving $\boldsymbol{\ell}$ in the gradient direction at each step along with projecting it onto $\mathcal{L}$, and the details are listed in Algorithm 1. Besides, in practice, the real load samples may be out of distribution, and incur  costs which are over the upper bound and thus causing screening failure. Therefore, we use a relaxation parameter $\epsilon$ to adjust the obtained upper bound $\texttt{PGA}(f_{NN}(\boldsymbol{\ell}))$ and integrate it to \eqref{screening1}. Then, we can get the following sample-agnostic screening problem considering the cost constraint,
\begin{subequations}  \label{Screening: sample-agnostic}
\begin{align}
\max_{\mathbf{u}, \mathbf{x}, \mathbf{f}, \boldsymbol{\ell}} / \min _{\mathbf{u}, \mathbf{x}, \mathbf{f},  \boldsymbol{\ell}}\quad  & f_{j} \label{Screening2:obj}\\
\text { s.t. } \quad & 
\eqref{Screening:gen} \eqref{Screening:flow}\eqref{Screening:balance}\eqref{Screening:u}\eqref{Screening:load}\\
& \sum_{i=1}^n c_i x_i \leq \texttt{PGA}(f_{NN}(\boldsymbol{\ell}))(1+\epsilon). \label{Screening: cost bound}
\end{align}
\end{subequations}

\begin{figure}[b]
\removelatexerror
\vspace{-1em}
\label{Algorithm: PGA}
  \renewcommand{\algorithmicrequire}{\textbf{Input:}}
  \renewcommand{\algorithmicensure}{\textbf{Output:}}
  \begin{algorithm}[H]
    \caption{Projected Gradient Ascent Algorithm}
    \begin{algorithmic}[1]
      \REQUIRE Load distribution $\mathcal{L}$, trained NN model $f_{NN}(\boldsymbol{\ell})$, step size $\beta$.
      \ENSURE Upper bound $\texttt{PGA}(f_{NN}(\boldsymbol{\ell}))$.
      \renewcommand{\algorithmicensure}{\textbf{Initialize:}}
      \ENSURE Random load vector $\boldsymbol{\ell}^{(0)} \in \mathcal{L}$, $k=0$.
      \WHILE {$\boldsymbol{\ell}^{(k)}$ doesn't converge}
      \STATE Update: $\texttt{PGA}(f_{NN}(\boldsymbol{\ell})) \leftarrow f_{NN}(\boldsymbol{\ell}^{(k)})$
      \STATE Calculate gradient $\nabla_{\boldsymbol{\ell}} f_{NN}(\boldsymbol{\ell})$
      \STATE Update: $\boldsymbol{\ell}^{(k+1)} \leftarrow \texttt{Proj}_{\mathcal{L}}(\boldsymbol{\ell}^{(k)}+\beta\nabla_{\boldsymbol{\ell}} f_{NN}(\boldsymbol{\ell}))$
      \STATE $k \leftarrow k+1$
      \ENDWHILE
      \STATE Return $\texttt{PGA}(f_{NN}(\boldsymbol{\ell}))$
    \end{algorithmic}
  \end{algorithm}
\end{figure}

In sample-aware case, we predict and still relax the UC cost for each specific sample, as the predicted cost may be lower than the actual cost, which can result in an infeasible adjusted screening problem for the investigated sample. Then we add the relaxed cost to \eqref{screening2}, and the adjusted sample-agnostic screening problem can be formulated as follows,
\begin{subequations} \label{Screening2: sample-aware}
\begin{align}
\max_{\mathbf{u}, \mathbf{x}, \mathbf{f}} / \min _{\mathbf{u}, \mathbf{x}, \mathbf{f}}\quad  & f_{j}\\
\text { s.t. } \quad & 
\eqref{Screening2:gen} \eqref{Screening2:flow}\eqref{Screening2:balance}\eqref{Screening2:u}\\
& \sum_{i=1}^n c_i x_i \leq f_{NN}(\boldsymbol{\ell})(1+\epsilon). \label{Screening2: cost bound}
\end{align}
\end{subequations}

Note that the upper bound given by the NN model will be a constant given load region or specific load vector, so the screening problems \eqref{Screening: sample-agnostic} and \eqref{Screening2: sample-aware} can be treated as linear programming problems which are efficient to solve.

\vspace{0.5em}
\section{Case Study}

\begin{figure*}[ht]
\vspace{-1.5em}
\centering
\includegraphics[width=0.9\linewidth]{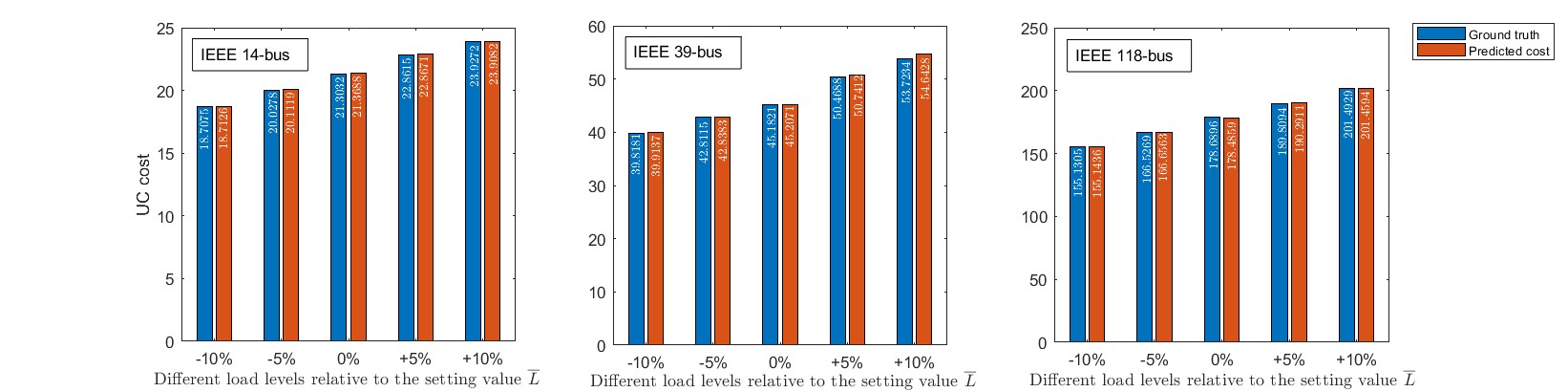}
\caption{The NN’s predicted costs of different load levels.} \label{Results: Prediction}
\vspace{-0.5em}
\end{figure*}

\begin{figure*}[ht]
\centering
\includegraphics[width=1.1\linewidth]{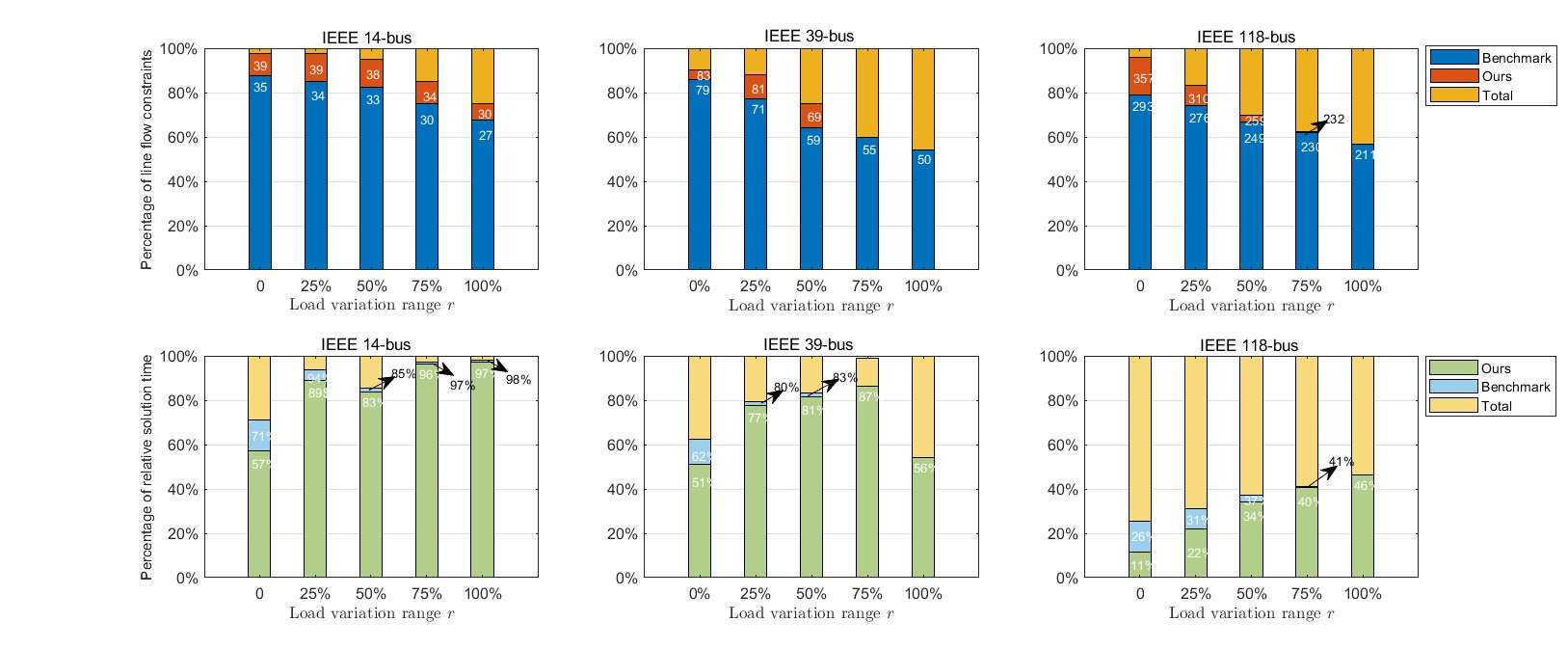}
\vspace{-2em}
\caption{Percentage of the reduced constraints (upper) and relative solution time (lower) of sample-agnostic screening on different load variation ranges.} \label{Results: Reduced_cons_solu_time_agnostic}
\vspace{-1em}
\end{figure*}

To evaluate the performance of the proposed constraint screening algorithm, we take the original optimization-based method as benchmark, and also compare our method with KNN method in this section. The predicting accuracy of the learning cost predictors, the computational efficiency and the solution accuracy of the reduced UC problem are examined over a wide range of problem settings. The details are given at \url{https://github.com/Hexuan085/UC_SCREENING_ML}.


\subsection{Simulation Setups}
We carry out the numerical experiments on IEEE 14-bus, IEEE 39-bus and IEEE 118-bus power systems. For each system, we consider the load with 0\%, 25\%, 50\%, 75\% and 100\% variation which is defined as $r$ around the average nominal values $\overline{\boldsymbol{\ell}}$. When investigating the sample-aware constraint screening, the load level is known in our setting and is defined as $\overline{L}$. Then the load region $\mathcal{L}$ considered here can be represented as:
\begin{subequations}
\begin{align}
(1-r)\overline{\boldsymbol{\ell}}\leq&\boldsymbol{\ell}\leq(1+r)\overline{\boldsymbol{\ell}}\label{UC:load range}\\
\sum_{i=1}^n&l_i=\overline{L}.\label{UC:load level}
\end{align}
\end{subequations}

\begin{table}[htbp]
\centering
\caption{Comparisons of the relative cost error and relative solution time for IEEE 118-bus system} \label{Results: Comparison__agnostic}
\setlength{\tabcolsep}{1.8mm}{\begin{tabular}{c|cccc|cccc}
\hline
      & \multicolumn{4}{c|}{Total cost error (\%)} & \multicolumn{4}{c}{Total solution time (\%)} \\ \hline
\diagbox{Method}{Range}& 25      & 50     & 75      & 100    & 25           & 50           & 75           & 100         \\ \hline
KNN5  & 5.3     & 9.8    & 3.3    & 4.5    & 18.3         & 16.5         & 16.4         & 16.7        \\
KNN10 & 0.8     & 0      & 0.2     & 1.7    & 18.5         & 19.1         & 19.7         & 20.3        \\
Benchmark    & 0       & 0      & 0       & 0      & 31.4         & 36.1         & 40.8         & 45.6        \\
Ours   & 0       & 0      & 0       & 0      & 21.7         & 33.4         & 39.5         & 45.7        \\ \hline
\end{tabular}}
\vspace{-1em}
\end{table}

\begin{table}[htbp] 
\caption{Percentage of average reduced constraints of Sample-aware screening} \label{Results: Reduced_constraints__aware}
\setlength{\tabcolsep}{1.6mm}{
\begin{tabular}{ccccccc}
\hline
        & Num.Gen. & Num.Lines & Benchmark & Ours & Actual& $\epsilon$ \\ \hline
14-bus  & 5        & 20        & 92.5\%     & 97.5\%        & 97.5\% & 0.01\\
39-bus  & 10       & 46        & 84.7\%     & 86.9\%       & 92.4\% & 0.03\\
118-bus & 54       & 186       & 81.5\%     & 85.7\%       & 97.3\% & 0.01\\ \hline
\vspace{-3em}
\end{tabular}}
\end{table}

To generate samples for training and validating the neural network model and KNN model, we use uniform distribution to get different load vectors $\boldsymbol{\ell} \in \mathcal{L}$ for sample-agnostic case or random $\boldsymbol{\ell}$ for sample-aware case, and then solve (\ref{UC}) for all loads. The UC cost and the binding situation of each line flow constraint are recorded. Under each setting, we solve and collect 10,000 samples for each neural network with 20 percentage of generated data split as test samples, while for KNN we solve 2,000 samples only based on 118-bus system due to computation burden. Moreover, when evaluating the screening performance of the benchmark, the proposed method and KNN, we use the same validation data and consider 100 samples for each validation case. The used neural networks all have ReLU activation units and 4 layers, and corresponding neurons on each hidden layer are 50, 30, 30. We feed the load vector as input for the neural network and the output is the corresponding UC cost, then we further use the cost to solve (\ref{Screening: sample-agnostic}) and (\ref{Screening2: sample-aware}).

All simulations have been carried out on a laptop with a 2.50 GHz processor and 16G RAM. Specifically, all the optimization problems are modeled using Python and solved with CVXPY\cite{Steven2016CVXPY} powered by GPLK\_MI solver \cite{makhorin2008glpk}.
\vspace{-0.1em}
\subsection{Simulation results}
To ensure the effectiveness and scalability of the proposed cost predictors, we train the NNs for different load levels, and randomly select a specific load vector from each load level to predict the corresponding cost. The results are shown in Fig. \ref{Results: Prediction}, from which it can be seen that the predicted costs are almost equal to the actual costs obtained by solving (\ref{UC}) with the relative error less than 1\%. Note that the predicted cost can be lower than the actual cost, so it is reasonable to consider $\epsilon$ in (\ref{Screening2: cost bound}) to ensure feasibility.

Using the NN models trained for the setting load levels and the PGA algorithm, we can get the upper bounds in (\ref{Screening: cost bound}) so as to conduct the sample-agnostic constraint screening. Then, this method is compared with the benchmark and KNN method, which is carried out on 118-bus system and the results are shown in Table \ref{Results: Comparison__agnostic} and Fig. \ref{Results: Reduced_cons_solu_time_agnostic}. According to Table \ref{Results: Comparison__agnostic} where the cost error and the solution time of the reduced problems are relative to the result of the original problem (\ref{UC}), KNN methods can reduce more solution time than other methods. The total cost errors of the case K=5 are lower than that of the case K=10, while the situation of the solution time is on the contrary. Though requiring more solution time, the benchmark and our methods can promise the solution accuracy without cost error. 

Meanwhile, our method can screen more constraints and save more solution time than the benchmark in all cases of investigated power systems with load variation range from 0\% to 50\% according to Fig. \ref{Results: Reduced_cons_solu_time_agnostic}. In the cases of 39-bus and 118-bus systems with 75\% to 100\% load range, the performance of the two methods are very close. This may be due to the increasing patterns of non-redundant constraints with larger load variation range, i.e., the percentage of the redundant constraints decreases when widening the load variation as shown in Fig. \ref{Results: Reduced_cons_solu_time_agnostic}.   
 
Furthermore, according to Table \ref{Results: Reduced_constraints__aware}, the average percentage of redundant constraints can reach 92.4\% to 97.3\% for a specific load vector. Then, with the sample-aware constraint screening, most of the redundant constraints can be removed. Specifically, the benchmark method can screen 81.5\% to 92.5\% of total constraints as redundant, while our method defined in (\ref{Screening2: sample-aware}) can screen 85.7\% to 95.1\% with setting $\epsilon$ properly. 

The above results show the following positive effects of our method:
\begin{enumerate}
\item Capturing the mapping between load vector and UC cost well at different load levels.
\item Realizing the trade-off between computational efficiency and solution accuracy in the sample-agnostic case.
\item Achieving higher screening efficiency in the sample-aware case.
\end{enumerate}
\vspace{0.5em}


\vspace{-0.8em}
\section{Conclusion and Future Works}
In this paper, we introduce a novel usage of machine learning to help screen redundant constraints. The neural networks are trained to predict UC cost so as to integrate the cost constraints to original screening problem efficiently. With the cost constraints, the search space of constraint screening can be sufficiently tightened. Since our method does not necessarily yield a minimal set of active constraints for the underlying UC problem, in the future work we would like to  seek theoretical understandings about the set of constraints and investigate how the proposed techniques can be generalized to multi-step UC problem with nonlinear constraints. We also plan to explore the potential of making sample-agnostic case serve as the warm-start for the sample-aware case.
\bibliographystyle{IEEEtran}
\bibliography{bib}
\end{document}